\documentclass[11pt]{article}
\usepackage{amssymb}
\usepackage{latexsym}

\def\QED{\hfill$\Box$}
\def\demo{\noindent{\bf Proof. }}

\newtheorem{Theorem}{Theorem}[section]
\newtheorem{Definition}[Theorem]{Definition}
\newtheorem{Lemma}[Theorem]{Lemma}
\newtheorem{Corollary}[Theorem]{Corollary}

\newtheorem{Remark}[Theorem]{Remark}
\newtheorem{Example}[Theorem]{Example\/}

\newtheorem{Proposition}[Theorem]{Proposition}

\begin{document}
\topmargin3mm
\medskip

\begin{center}
{\large\bf Explicit representations by halfspaces of the \\ edge cone of a graph}
\vspace{6mm}\\
{Carlos E. Valencia
 \ and Rafael H. Villarreal\footnote{Partially supported by SNI, M\'exico.}}


{\small Departamento de Matem\'aticas}\\ 
{\small Centro de Investigaci\'on y de Estudios Avanzados del
IPN}\\   
{\small Apartado Postal 14--740}\\ 
{\small 07000 M\'exico City, D.F.}\\ 
{\small e-mail {\tt vila@esfm.ipn.mx}}\\
\end{center}
\date{}


\medskip

\begin{abstract}
\noindent Let $G$ be an arbitrary graph. The main results are explicit 
representations of the edge cone 
of $G$ as a finite intersection of closed halfspaces. If $G$ is bipartite 
and connected we determine the facets of the edge cone and present a 
canonical irreducible representation. 
\end{abstract}

\vspace{0.5cm}

\section{Introduction}

Let $G$ be an arbitrary graph on the vertex set $V=\{v_1,\ldots,v_n\}$. 
The {\em edge cone\/} of 
$G$ is the cone ${\mathbb R}_+{\cal A}
\subset {\mathbb R}^n$ spanned by
the set ${\cal A}$ of all vectors $e_i+e_j$ such that $v_i$ is
adjacent to $v_j$, where $e_i$ denotes the 
$i${\it th} unit vector.  

Our first main goal is to give an explicit combinatorial 
description of the edge cone of $G$, see 
Theorem~\ref{feasible1} and Corollary~\ref{feasible2}. 
This description generalizes that of \cite[Corollary~3.3]{join}. 
In loc. cit. only the non bipartite case 
was studied. 

The second main goal is to study in detail the facets of the edge 
cone of a connected bipartite graph and show a canonical irreducible representation 
of the edge cone, see Proposition~\ref{dec27-01} and Theorem~\ref{uniqrepf}. 
As an application the classical marriage theorem will follow. Our results can 
be applied to commutative 
algebra to compute the $a$-invariant 
and the canonical module of an edge subring along the lines of
\cite{BVV,git-val,shiftcon,monalg}. It is a bit 
surprising that the edge cone of a bipartite graph has not been studied 
before from the point of view of polyhedral geometry. 

To show our results we use graph theory,  linear algebra (Farkas's Lemma, incidence 
matrices of graphs, Carath\'eodory's Theorem, Kronecker's Lemma), and 
polyhedral geometry (finite basis theorem and facet structure). The proofs 
require a careful analysis at the graph theoretical level. Our main 
references for graphs, algebra and geometry are \cite{Boll,Bron,ewald,Har,Schr}.

\section{Preliminaries} 



Let ${\cal A}=\{\alpha_1,\ldots,\alpha_q\}$ be a finite set of vectors 
in $\mathbb{Q}^n$. The {\em cone\/} $\mathbb{R}_+{\cal A}$ generated 
by the set ${\cal A}$ is defined as
\[
{\mathbb R}_+{\cal A}:=
\left.\left\{\sum_{i=1}^qa_i\alpha_i\right\vert\, a_i\in
{\mathbb R}_+\, \mbox{ for all }i\right\}\subset\mathbb{R}^n,
\]
where $\mathbb{R}_+$ is the set of non negative real numbers.

By the finite basis theorem \cite[Chapter~4]{webster} ${\mathbb R}_+{\cal A}$ 
is a rational polyhedral cone, that is, ${\mathbb R}_+{\cal A}$  is the intersection 
of finitely many {\it closed halfspaces\/} of the form:
$$
H_a^+:=
\{x\in\mathbb{R}^n\vert\, \langle x,a\rangle\geq 0\},
$$
where $0\neq a\in\mathbb{Z}^n$ and the nonzero entries of $a$ are relatively 
prime. Here $\langle x,a\rangle$ denotes the standard 
inner product of 
$x$ and $a$. 

Note that if $H_a^-:=H_{-a}^+$, then the intersection $H_a^+\cap H_a^-$ is the 
bounding hyperplane 
$$
H_a:=\{x\, \vert\, \langle x,a\rangle=0\}
$$ 
with normal vector $a$.

To simplify notation set $Q=\mathbb{R}_+{\cal A}$. Recall that a subset $F\subset\mathbb{R}^n$ 
is a {\it face\/} of $Q$ if $F=Q\cap H_a$ for some hyperplane $H_a$ such that 
$Q\subset H^+_a$ or $Q\subset H_a^-$. The 
hyperplane $H_a$ is called a {\it supporting hyperplane\/} of $Q$. The 
{\it improper faces\/} of $Q$ are $Q$ and $\emptyset$, all the other faces 
are called {\it proper faces\/}. If a face of $Q$ has dimension 
 $\dim(Q)-1$ it is called a {\em facet\/}. The dimension of $Q$ 
is by definition the dimension of ${\rm aff}(Q)$, the affine hull of $Q$. Note 
that a face of $Q$ is again a finitely generated cone, see \cite{Bron}. 

\begin{Definition}\rm If a polyhedral cone $Q=\mathbb{R}_+{\cal A}$ in $\mathbb{R}^n$ is 
represented as 
$$Q={\rm aff}(Q)\cap\left(\bigcap_{i=1}^rH^+_{a_i}\right)\eqno(*)$$
and satisfies  
$$
Q\neq {\rm aff}(Q)\cap\left(\bigcap_{i\neq j}^rH^+_{a_i}\right)
$$
for all $j$, we say that $(*)$ is an 
{\it irreducible representation\/}
\index{irreducible!representation} of $Q$. 
\end{Definition}


Part of the importance of an irreducible representation can be seen in the 
following general fact.

\begin{Theorem}\label{rep-of-facets} Let $Q$ be a polyhedral 
cone in 
${\mathbb R}^n$ which is not an affine space. If      
$$
Q={\rm aff}(Q)\cap 
H^+_{a_1}\cap\cdots\cap H^+_{a_r}
$$
is an irreducible representation of $Q$ with $a_i\in\mathbb{R}^n\setminus\{0\}$ 
for all $i$, then the facets of $Q$ are precisely 
the sets $F_1,\ldots,F_r$, where $F_i=Q\cap H_{a_i}$. Moreover 
each proper face of $Q$ is the intersection of those 
facets of $Q$ that contain it.
\end{Theorem}

\demo See \cite[Theorem~3.2.1]{webster}. \QED 

\paragraph{The incidence matrix of a graph} 


In the sequel we use standard terminology and notation from 
graph theory and adopt the book of Harary \cite{Har} as our main reference. For 
the reader's convenience we recall a few notions about graphs.

\medskip

Let $G$ be a simple graph with vertex set 
$V(G)=\{v_1,\ldots,v_n\}$ and edge set 
$E(G)=\{z_1,\ldots,z_q\}$, thus every 
edge $z_i$ is an unordered pair of distinct vertices $z_i=\{v_{i_j},v_{i_k}\}$. 
The {\it incidence matrix\/} $M_G=[a_{ij}]$ associated to $G$ is the $n \times q$ matrix 
defined by
$$
a_{ij}=\left\{\begin{array}{ll}
1&\mbox{if }v_i\in z_j, \mbox{ and}\\
0&\mbox{if }v_i\notin z_j.
\end{array} \right. 
$$
Note that each column of $M_G$ has exactly two $1$'s and 
the rest of its entries equal to zero. If 
$z_i=\{v_{i_j},v_{i_k}\}$ define $\alpha_i=e_{i_j}+e_{i_k}$, 
where $e_i$ is the $i${\it th} unit vector 
in $\mathbb{R}^n$. Thus the columns of $M_G$ are 
precisely the vectors $\alpha_1,\ldots,\alpha_q$. As an example consider a 
triangle $G$ with vertices 
$v_1,v_2,v_3$. In this case:
$$
M_G=\left(
\begin{array}{ccc}
1&0&1\\
1&1&0\\
0&1&1
\end{array}
\right),
$$
with the vectors $\alpha_1=e_1+e_2$, $\alpha_2=e_2+e_3$, and 
$\alpha_3=e_1+e_3$ corresponding to 
the edges $z_1=\{v_1,v_2\}$, $z_2=\{v_2,v_3\}$, and 
$z_3=\{v_1,v_3\}$.

Recall that a graph $G$ is {\it bipartite\/} if there is a {\it bipartition\/} 
$(V_1,V_2)$ of $G$, that is, $V_1$ and $V_2$ are vertex classes 
satisfying:
\begin{description}
\item{\rm (a)} $V(G)=V_1\cup V_2$,\vspace{-1mm}
\item{\rm (b)} $V_1\cap V_2=\emptyset$, and\vspace{-1mm}
\item{\rm (c)} every edge of $G$ joins a vertex of $V_1$ to a vertex of $V_2$.
\end{description}
If $G$ is connected such a bipartition is uniquely determined. Equivalently $G$ 
is bipartite if all its cycles are of even length.  If $G$ is bipartite, then 
its incidence matrix is totally unimodular, that is, all the $i\times i$ 
minors of $M_G$ are equal to $0$ or $\pm 1$ for all $i\geq 1$, see \cite{Schr}. 

\section{An explicit representation of the edge cone} 

Let us introduce some 
more terminology and fix some more notation that
will be used throughout. 

Let $G$ be a simple graph and let $M_G$ be its incidence matrix. 
We set ${\cal A}_G$ (or simply $\cal A$ if $G$ is understood) equal
to the set $\{\alpha_1,\ldots,\alpha_q\}$ of column vectors  of $M_G$. Since 
$\alpha_i$ represents an edge of $G$ sometimes $\alpha_i$ is called an edge 
or an {\it edge vector\/}. The {\em edge cone\/}  of $G$ is defined as the cone 
${\mathbb R}_+{\cal A}$ generated by ${\cal A}$. Note ${\mathbb R}_+{\cal A}\neq(0)$ if 
$G$ is not a discrete graph. By \cite{Kulk1} one has 
$$
n-c_0(G)={\rm rank}(M_G)=\dim\, {\mathbb R}_+{\cal A},
$$
where $c_0(G)$ is the number of bipartite connected components of $G$. 

\begin{Lemma}
If $v_i$ is not an isolated vertex of $G$, then the 
set $F=H_{e_i}\cap {\mathbb R}_+{\cal A}$ is a proper 
face of the edge cone.
\end{Lemma} 

\demo Note $F\neq\emptyset$ because $0\in F$, and 
${\mathbb R}_+{\cal A}\subset H_{e_i}^+$. Since $v_i$ is not 
an isolated vertex ${\mathbb R}_+{\cal A}\not\subset H_{e_i}$. \QED

\medskip

Given a subset $A\subset {V}(G)$, the {\it neighbor set\/} of $A$, denoted 
$N_G(A)$ or simply $N(A)$, is defined as
$$
N(A)=\{v\in {V}(G)\mid v \; \mbox{\rm is adjacent to some vertex
in}\; A\}.
$$ 

Let $A$ be an {\it independent set\/} of vertices of $G$, that is, no two vertices 
of $A$ are adjacent. The supporting
hyperplane of the edge cone of $G$ defined by  
\[
\sum_{v_i\in A}x_i=\sum_{v_i\in N(A)}x_i
\]
will be denoted by $H_A$. 


\begin{Lemma} If $A$ is an independent set of 
vertices of $G$ and $F={\mathbb R}_+{\cal A}\cap H_A$, 
then either $F$ is a proper face of the edge cone or 
$F={\mathbb R}_+{\cal A}$.
\end{Lemma}

\demo It suffices to prove the containment 
$\mathbb{R}_+{\cal A}\subset H_A^-$. Take and 
edge $\{v_j,v_\ell\}$ of $G$. 
If $\{v_j,v_\ell\}\cap A\neq\emptyset$, then $e_j+e_\ell$ 
is in $H_A$, else $e_j+e_\ell$ is in $H_A^-$. \QED

\begin{Definition}\rm The {\it support} of a vector 
$\beta=(\beta_i)\in{\mathbb R}^n$ is defined as 
$${\rm
supp}(\beta)=\{\beta_i\, |\, \beta_i\neq 0\}.$$ 
\end{Definition}

\begin{Lemma}[\cite{join}]\label{principal}
Let $V=\{v_1,\ldots,v_n\}$ be the vertex set of $G$ and let $G_1,\ldots,G_r$ be the 
 connected components of $G$. If 
$G_1$ is a tree with at least two vertices 
and $G_2,\ldots,G_r$ are unicyclic non bipartite graphs, then 
${\rm ker}(M_G^t)=(\beta)$ for some $\beta$ in $\mathbb{R}^n$ with 
${\rm supp}(\beta)=\{1,-1\}$ such that 
$V(G_1)=\{v_i\in V\, \vert\, \beta_i=\pm 1\}$. 
\end{Lemma}
For use below we recall the following form of Farkas's Lemma, which is 
called the fundamental theorem of linear inequalities, see \cite[Theorem~7.1]{Schr}. 
\begin{Theorem}\label{ftli} Let ${\cal A}=
\{\alpha_1,\ldots,\alpha_q\}$ be a set of vectors in $\mathbb{R}^n$ 
and let $\alpha\in\mathbb{R}^n$. If 
$\alpha\notin\mathbb{R}_+{\cal A}$ and 
$t={\rm rank}\{\alpha_1,\ldots,\alpha_q,\alpha\}$, then there exists a hyperplane 
$H_a$ containing $t-1$ linearly independent vectors from ${\cal A}$ 
such that $\langle a,\alpha\rangle>0$ and $\langle a,\alpha_i\rangle\leq 0$ 
for $i=1,\ldots,q$.
\end{Theorem}

\begin{Theorem}\label{feasible1}
If $G$ is a connected graph with vertex set
$V=\{v_1,\ldots,v_n\}$ and $\mathbb{R}_+{\cal A}$ 
is the edge cone of $G$, then 
$$
{\mathbb R}_+{\cal A}=\left(\bigcap_{A\in\cal F} H_A^-\right)
\bigcap\left(\bigcap_{i=1}^n H_{e_i}^+\right),\eqno(*)
$$
where $\cal F$ is the family of all the independent sets 
of vertices of $G$ and $H_{e_i}^+$ is the closed 
halfspace $\{x\in\mathbb{R}^n\, \vert\, x_i\geq 0\}$.  
\end{Theorem}
\demo Let ${\cal A}=\{\alpha_1,\ldots,\alpha_q\}$ be the set of
column vectors of the incidence matrix of $G$. Since $\mathbb{R}_+{\cal A}$ is 
clearly contained in the right hand 
side of Eq.~$(*)$ it suffices to prove the other containment. Take 
$\alpha\in\mathbb{R}^n$ in the right hand side of Eq.~$(*)$. The proof 
is by contradiction, that is, assume $\alpha\notin {\mathbb R}_+{\cal A}$. 
By \cite[Corollary~3.3]{join} we may assume $G$ bipartite with $n\geq 3$ 
vertices. 

Note that if $(V_1,V_2)$ is the bipartition of $G$, 
then ${\rm aff}({\mathbb R}_+{\cal A})$ is the hyperplane
$$
\sum_{v_i\in V_1}x_i=\sum_{v_i\in V_2}x_i,
$$
because ${\rm dim}(\mathbb{R}_+{\cal A})=n-1$. As $H_{V_1}^-\cap H_{V_2}^-=H_{V_1}$, 
the vector $\alpha$ is in 
${\rm aff}({\mathbb R}_+{\cal A})$. As a consequence 
${\rm rank}({\cal A}\cup\{\alpha\})=n-1$. 

By Theorem~\ref{ftli} there is $a\in\mathbb{R}^n$ and there are linearly 
independent  vectors $\alpha_1,\ldots,\alpha_{n-2}$ 
in ${\cal A}$ such that \vspace{-1mm}
\begin{description}
\item{\rm (i)\ \,} 
$\langle a,\alpha_i\rangle=0 \, \,  \mbox{ for } i=1,2,\ldots,n-2$,
\vspace{-1mm}
\item{\rm (ii)\,} $\langle a,\alpha_i\rangle\leq 0
\, \,\mbox{ for } i=1,2,\ldots,q$, and \vspace{-1mm}
\item{\rm (iii)} $\langle a,\alpha\rangle>0$.
\end{description}

Observe that $\mathbb{R}_+{\cal A}\not\subset H_a$ because 
${\rm aff}(\mathbb{R}_+{\cal A})\neq H_a$. There exists $\alpha_j$ in 
${\cal A}$ such that $\alpha_1,\ldots,\alpha_{n-2},\alpha_j$ is a basis 
of ${\rm aff}(\mathbb{R}_+{\cal A})$ 
as a real vector space. In particular we can write 
\begin{equation}\label{jan27-02}
\alpha=\lambda_1\alpha_1+\cdots+\lambda_{n-2}\alpha_{n-2}+
\lambda_j\alpha_j\ \ \ (\lambda_i\in\mathbb{R})
\end{equation}
It follows that $\langle \alpha,a\rangle=
\lambda_j\langle \alpha_j,a\rangle>0$. Thus $\lambda_j<0$. 

Consider the subgraph $D$ of $G$ whose
edges correspond to  
$\alpha_1,\ldots,\alpha_{n-2}$ and its vertex set is 
the union of the vertices in the edges of $D$. Set $k=|V(D)|$. By
 \cite{Kulk1} one has:
$$
n-2={\rm rank}(M_D)=k-c_0(D),
$$
where $M_D$ is the incidence matrix of $D$ and $c_0(D)$ is the number of 
bipartite components of $D$. Thus $0\leq n-k=2-c_0(D)$. This shows that 
either $c_0(D)=1$ and $k=n-1$ or $c_0(D)=2$ and $k=n$.

Case (I): Assume $C_0(D)=1$ and $k=n-1$. Set $V(D)=\{v_1,\ldots,v_{n-1}\}$. 
As $D$ is a tree with $n-2$ edges and $\langle\alpha_i,a\rangle=0$ 
for $i=1,\ldots,n-2$, applying Lemma~\ref{principal}, 
one may assume $a=(a_1,\ldots,a_{n-1},a_n)$, where $a_i=\pm 1$ 
for $1\leq i\leq n-1$. 

Set $a'=(0,\ldots,0,-1)=-e_n$. Next we prove the following 
\begin{description}
\item{\rm (a)} $\langle \alpha_i,a'\rangle=0$ for $i=1,\ldots,n-2$.\vspace{-1mm}
\item{\rm (b)} $\langle \alpha_j,a'\rangle=-1$ and $\langle \alpha_j,a\rangle<0$. 
\end{description}
Condition (a) is clear. To prove (b) first note 
$\alpha_j\notin\mathbb{R}(\alpha_1,\ldots,\alpha_{n-2})$. Then 
$\alpha_j=e_k+e_n$, because otherwise the ``edge'' $\alpha_j$ added to 
the tree $D$ form a graph with a unique even cycle, to derive a contradiction 
recall that a set of ``edge vectors'' forming an even cycle  are linearly 
dependant.  Thus $\langle \alpha_j,a'\rangle=-1$. 

On the other hand 
$\langle \alpha_j,a\rangle<0$, because if $\langle \alpha_j,a\rangle=0$, 
then the hyperplane $H_a$ would contain the linearly independent vectors 
$\alpha_1,\ldots,\alpha_{n-2},\alpha_j$ and consequently 
${\rm aff}(\mathbb{R}_+{\cal A})$ would be equal to $H_a$, a contradiction. 

To finish the proof of this case we use the inequality
$$
\langle\alpha,a'\rangle=\lambda_j\langle\alpha_j,a'\rangle>0
$$
to conclude $\langle\alpha,a'\rangle>0$, a contradiction because 
$\alpha\in H_{e_n}^+$. 

Case (II): Assume $C_0(D)=2$ and $k=n$. Let $D_1$ and $D_2$ be the 
components of $D$ and set $U_1=V(D_1)$ and $U_2=V(D_2)$.  

Using Lemma~\ref{principal} we can relabel the vertices of the graph $D$ and 
write $a=rb+sc$, where $0\neq r\geq s\geq 0$ are rational numbers,
$$
b=(b_1,\ldots,b_m,0,\ldots,0),\ \ \ \ c=(0,\ldots,0,c_{m+1},\ldots,c_n),
$$
$U_1=\{v_1,\ldots,v_m\}$, $b_i=\pm 1$ for $i\leq m$, 
and $c_i=\pm 1$ for $i>m$. Set $a'=b$. Note the following:
\begin{description}
\item{\rm (a)} $\langle \alpha_i,a'\rangle=0$ for $i=1,\ldots,n-2$.\vspace{-1mm}
\item{\rm (b)} $\langle \alpha_j,a'\rangle=-1$ and $\langle \alpha_j,a\rangle<0$; 
this holds for any $\alpha_j\notin\mathbb{R}(\alpha_1,\ldots,\alpha_{n-2})$.
\end{description}
Condition (a) is clear. To prove (b) first note that the 
inequality $\langle \alpha_j,a\rangle<0$ can be shown as in case (I). Observe that 
if an ``edge'' $\alpha_k$ has vertices in $U_1$ (resp. $U_2$), then 
$\langle \alpha_k,a\rangle=0$. Indeed if we add the edge $\alpha_k$ to the tree $D_1$ 
(resp. $D_2$) we get a graph with a unique even cycle and this implies 
that $\alpha_1,\ldots,\alpha_{n-2},\alpha_k$ are linearly dependant, that is, 
$\langle \alpha_k,a\rangle=0$. Thus $\alpha_j=e_i+e_\ell$ for some $v_i\in U_1$ 
and $v_\ell\in U_2$. From the inequality
$$
\langle \alpha_j,a\rangle=r\langle \alpha_j,b\rangle+s\langle \alpha_j,c\rangle=
rb_i+sc_\ell<0
$$   
we obtain $b_i=-1=\langle\alpha_j,a'\rangle$, as required.   

Next we set 
\[
A=\{v_i\in V\vert\, b_i=1\}\mbox{\ \ and\ \ }B=\{v_i\in V\vert\,
b_i=-1\}.
\]
Note that $\emptyset\neq A\subset U_1$ and $\emptyset\neq B\subset U_1$, 
because $D_1,D_2$ are trees with at least two vertices. 
We will show that $A$ is an independent set of 
$G$ and $B=N_G(A)$. 

If $A$ is not an independent set of $G$, there is an edge $\{v_i,v_\ell\}$ of $G$ 
for some $v_i,v_\ell$
in $A$. Thus $\alpha_k=e_i+e_\ell$, by (a) and (b) we get 
$\langle a',\alpha_k\rangle\leq 0$, which is impossible because $\langle
a',\alpha_k\rangle=2$. This proves that $A$ is an independent set of $G$.

Next 
we show $N_G(A)=B$. If $v_i\in N_G(A)$, then $\alpha_k=e_i+e_\ell$ for
some $v_\ell$ in $A$, using (a) and (b) we obtain $\langle 
a',\alpha_k\rangle=b_i+1\leq 0$ and $b_i=-1$, hence $v_i\in B$.
Conversely if $v_i\in B$, since $D_1$ has no isolated vertices, there
is $1\leq k\leq n-2$ so that $\alpha_k=e_i+e_\ell$, for some 
$\ell$, by (b) we obtain $\langle a',\alpha_k\rangle=-1+b_\ell=0$,
which shows that $v_\ell\in A$ and $v_i\in N_G(A)$.  

Therefore $H_{a'}=H_A$. Since 
$\mathbb{R}_+{\cal A}\not\subset H_{a'}$ (this follows from (b)), 
there is $\alpha_\ell\notin H_A$, thus $H_{a'}^-\cap H_A^-\neq\emptyset$ 
and consequently $H_A^-=H_{a'}^-$. By hypothesis $\alpha\in H_A^-$, 
hence $\langle\alpha,a'\rangle\leq 0$. From Eq.~(\ref{jan27-02}) together 
with and (a) and (b) one has
$$
\langle\alpha,a'\rangle=\lambda_j\langle\alpha_j,a'\rangle=-\lambda_j>0,
$$
a contradiction.  \QED 

\medskip

The next two results give an explicit representation by closed halfspaces of the edge cone of an 
arbitrary graph. Those representations were known for connected non 
bipartite graphs only \cite{join}. 

\begin{Corollary}\label{oct12-01} If $G$ is a graph with vertex set
$V=\{v_1,\ldots,v_n\}$ and $\mathbb{R}_+{\cal A}$ 
is the edge cone of $G$, then 
$$
{\mathbb R}_+{\cal A}=\left(\bigcap_A H_A^-\right)
\bigcap\left(\bigcap_{i=1}^n H_{e_i}^+\right),\eqno(*)
$$
where the intersection is taken over all the independent sets 
of vertices $A$ of $G$ and $H_{e_i}^+=\{x\in\mathbb{R}^n\, \vert\, x_i\geq 0\}$.  
\end{Corollary}

\demo Let $G_1,\ldots,G_r$ be the connected components of $G$. For 
simplicity of notation we assume $r=2$ and $V(G_1)=\{v_1,\ldots,v_m\}$. 
There is a decomposition
$$
{\mathbb R}_+{\cal A}=
{\mathbb R}_+{\cal A}_{G_1}\oplus{\mathbb R}_+{\cal A}_{G_2}.
$$
Let $\delta$ be a vector in the right hand side of Eq.~$(*)$. One can write
$$
\delta=(\delta_i)=\beta+\gamma=(\delta_1,\ldots,\delta_m,0,\ldots,0)+
(0,\ldots,0,\delta_{m+1},\ldots,\delta_n).
$$
Let $A$ be an independent set of $G_1$. Note $N_G(A)=N_{G_1}(A)$, hence 
$$
\sum_{v_i\in A}\delta_i\leq \sum_{v_i\in N_G(A)}\delta_i=
\sum_{v_i\in N_{G_1}(A)}\delta_i.
$$
Applying Theorem~\ref{feasible1} yields 
$\beta\in{\mathbb R}_+{\cal A}_{G_1}$. Similarly one has 
$\gamma\in{\mathbb R}_+{\cal A}_{G_2}$. Hence 
$\delta\in{\mathbb R}_+{\cal A}_{G}$, as required. \QED

\begin{Corollary}\label{feasible2}
Let $G$ be a graph with vertex set
$V=\{v_1,\ldots,v_n\}$. Then a vector $x=(x_1,\ldots,x_n)\in {\mathbb R}^n$ 
is in ${\mathbb R}_+{\cal A}$ if and only if $x$ is a solution of the 
system of linear inequalities 
\[
\begin{array}{rcll}
-x_i&\leq & 0, &i=1,\ldots,n
\\ & & &\vspace{-3mm} \\ 
\sum_{v_i\in A}x_i-\sum_{v_i\in N(A)}x_i& \leq & 0,& \mbox{for
all independent sets }A\subset V.
\end{array}
\] 
\end{Corollary}

\demo It follows at once from Corollary~\ref{oct12-01}. \QED

\medskip

\begin{Theorem}\label{oct13-01}
If $G$ is a graph with vertex set 
$V=\{v_1,\ldots,v_n\}$ and $F$ is a facet of the edge cone 
of $G$, then either 
\begin{description}
\item{\rm (a)} $F=\mathbb{R}_+{\cal A}\cap\{x\in\mathbb{R}^n\, \vert\, x_i=0\}$ for some 
$1\leq i\leq n$, or\vspace{-1mm} 

\item{\rm (b)} $F=\mathbb{R}_+{\cal A}\cap H_A$ for some independent set $A$ of $G$. 
\end{description}
\end{Theorem}

\demo By Corollary~\ref{oct12-01} we can write 
$${\mathbb R}_+{\cal A}={\rm aff}({\mathbb R}_+{\cal A})\cap 
H_1^-\cap\cdots\cap H_{r}^-$$
for some hyperplanes $H_1,\ldots,H_r$ such that none of the 
halfspaces $H_j^-$ can be omitted in the intersection and each $H_j$ is 
either of the form $H_{-e_i}$ or $H_j=H_A$ for some independent set $A$. 
By Theorem~\ref{rep-of-facets} the facets of ${\mathbb R}_+{\cal A}$ 
are precisely the sets $F_1,\ldots,F_r$, where 
$F_i=H_{i}\cap{\mathbb R}_+{\cal A}$. \QED

\begin{Remark}\rm (i) To verify whether a face $F$ as in (a) or (b) 
is a facet consider the set ${\cal B}$ of all $v_i\in\cal{A}$ 
that are in $F$. Note that $F$ is a facet if and only if 
$\dim\, \mathbb{R}_+{\cal B}=r-1$, where $r$ is the dimension 
of the edge cone. 

(ii) If $G$ is bipartite and connected, a graph theoretical characterization of 
the facets of the edge cone of $G$ will be given in Proposition~\ref{dec27-01}. 
The facets of the edge cone of $G$ for $G$ non bipartite were characterized 
in \cite[Theorem~3.2]{join}. 
\end{Remark}

\section{Studying the bipartite case}\label{jan30-02} For 
connected bipartite graphs we will present sharper 
results on the irreducible representations 
of edge cones and give a characterization of their facets.

\begin{Proposition}\label{oct20-01} Let $G$ be a connected bipartite graph 
with bipartition $(V_1,V_2)$. If $A$ is an independent set of $G$ such 
that $A\neq V_i$ for $i=1,2$, then $F={\mathbb R}_+{\cal A}\cap H_A$ is 
a proper face of the edge cone.
\end{Proposition}

\demo Assume $N(A)=V_2$. Take any $v_i\in V_1\setminus A$ and any 
$v_j\in V_2$ adjacent to $v_i$, then $e_i+e_j\notin H_A$. Thus 
we may assume $N(A)\neq V_i$ for $i=1,2$. 

Case (I): $N(A)\cap V_i\neq\emptyset$ for $i=1,2$. If the vertices in 
$N(A)\cap V_i$ for $i=1,2$ are only adjacent to vertices in $A$, then 
pick vertices $v_i\in N(A)\cap V_i$ and note that there is no path 
between $v_1$ and $v_2$, a contradiction. Thus there must be a vector 
in the edge cone which is not in $H_A$. 

Case (II): $A\subsetneq V_1$. If the vertices in $N(A)$ are only 
adjacent to vertices in $A$. Then a vertex in $A$ cannot 
be joined by a path to a vertex in $V_2\setminus N(A)$, a contradiction. As 
before we obtain $\mathbb{R}_+{\cal A}\not\subset H_A$. \QED

\begin{Proposition}\label{oct27-01} Let $G$ be a connected bipartite graph 
with bipartition $(V_1,V_2)$ and $\cal F$ the family of independent sets $A$ of $G$ such that 
$H_A\cap \mathbb{R}_+{\cal A}_G$ is a facet. If $A$ is in 
$\cal F$ 
and $V_i\cap A\neq \emptyset$ for $i=1,2$, then the halfspace $H_A^{-}$ 
is redundant in the following expression of the edge cone
$$
{\mathbb R}_+{\cal A}={\rm aff}(\mathbb{R}_+{\cal A})\cap
\left(\bigcap_{A\in\cal F} H_A^-\right)
\bigcap\left(\bigcap_{i=1}^n H_{e_i}^+\right).
$$
\end{Proposition}

\demo Set ${\cal A}=\{\alpha_1,\ldots,\alpha_q\}$. 
One can write $A=A_1\cup A_2$ with $A_i\subset V_i$ for $i=1,2$. There are 
$\alpha_1,\ldots,\alpha_{n-2}$ linearly independent vectors in 
$H_A\cap\mathbb{R}_+{\cal A}$, where $n$ is the number of vertices of $G$. 
Consider the subgraph $D$ of $G$ whose
edges correspond to  
$\alpha_1,\ldots,\alpha_{n-2}$ and its vertex set is 
the union of the vertices in those edges. Note that $D$ cannot 
be connected. Indeed there is no edge of $D$ connecting 
a vertex in $N_G(A_1)$ with a vertex in $N_G(A_2)$ because 
all the vectors $\alpha_1,\ldots,\alpha_{n-2}$ satisfy the equation
$$
\sum_{v_i\in A}{x_i}=\sum_{v_i\in N_G(A)}{x_i}.
$$
Hence by the proof of Theorem~\ref{feasible1} it follows that $D$ is a 
spanning subgraph of $G$ with 
two connected components $D_1$ and $D_2$ (which are trees) such that 
$V(D_i)=A_i\cup N_G(A_i)$, $i=1,2$. Therefore $H_{A_i}$ is a proper 
support hyperplane defining a facet $F_i=H_{A_i}\cap\mathbb{R}_+{\cal A}$, that 
is $A_1,A_2$ are in $\cal F$. Since $H_{A_1}^-\cap H_{A_2}^-$ is 
contained in $H_A^-$ the proof is complete. \QED

\begin{Proposition}\label{oct28-01} Let $G$ be a 
connected bipartite graph 
with bipartition $(V_1,V_2)$. If $A_2\subsetneq V_2$ and 
$F=H_{A_2}\cap\mathbb{R}_+{\cal A}$ is a facet of the edge cone of $G$, then 
$$
H_{A_2}^-\cap{\rm aff}({\cal A}')=\left\{\begin{array}{ll}
H_{A_1}^-\cap{\rm aff}({\cal A}')&\mbox{ where }A_1=V_1\setminus N(A_2)\neq\emptyset, \mbox{or}\\ 
H_{e_i}^+\cap{\rm aff}({\cal A}')&\mbox{for some vertex }
 v_i\mbox{ with }G\setminus\{v_i\}\mbox{ connected}, 
\end{array}
\right.
$$
where ${\cal A}'={\cal A}\cup\{0\}$.
\end{Proposition}

\demo Let us assume $G$ has $p$ vertices 
$v_1,\ldots,v_p$ and $V_1$ 
is the set of the first $m$ vertices of $G$. 
Set ${\cal A}=\{\alpha_1,\ldots,\alpha_q\}$. 
There are 
$\alpha_1,\ldots,\alpha_{p-2}$ linearly independent vectors in 
the hyperplane $H_{A_2}$. Consider the subgraph $D$ of $G$ whose
edges correspond to  
$\alpha_1,\ldots,\alpha_{p-2}$ and its vertex set is 
the union of the vertices in those edges. As $G$ is connected either $D$ is 
a tree with $p-1$ vertices or $D$ is a spanning subgraph of $G$ with two 
connected components.

If $D$ is a tree, write $V(D)=V(G)\setminus\{v_i\}$ for some $i$. Note 
$$
\langle \alpha_j,\alpha_{A_2}\rangle=-\langle \alpha_j,e_i\rangle\ \ \ \ 
(j=1,\ldots,q),
$$
where 
$$
\alpha_{A_2}=\sum_{v_i\in A_2}{e_i}-\sum_{v_i\in N(A_2)}{e_i}.
$$ 
Indeed if the ``edge'' $\alpha_j$ has vertices in $V(D)$, then both 
sides of the equality are zero, otherwise write $\alpha_j=e_i+e_\ell$. 
Observe $v_i\notin A_2$ and $v_\ell\in N(A_2)$ because $H_{A_2}$ being a facet cannot 
contain $\alpha_j$, thus both sides of the equality are equal to $-1$. As 
a consequence since ${\rm aff}({\cal A}')=
\mathbb{R}(\alpha_1,\ldots,\alpha_{p-2},\alpha_j)$ for some 
$\alpha_j=e_i+e_\ell$ we rapidly obtain
$$
\langle \alpha,\alpha_{A_2}\rangle=-\langle\alpha,e_i\rangle\ \ \ \ (\forall\, 
\alpha\in{\rm aff}({\cal A}')).
$$
Therefore
$$
H_{A_2}^-\cap{\rm aff}({\cal A}')=H_{e_i}^+\cap{\rm aff}({\cal A}'),
$$
as required. 

We may now assume $D$ is not a tree. We claim 
$A_1=V_1\setminus N(A_2)\neq\emptyset$. If $V_1=N(A_2)$. Take 
$v_i\in V_2\setminus A_2$ and $\{v_i,v_j\}$ and edge of $D$ containing 
$v_i$. Hence since $v_j\in N(A_2)$ we get 
$\langle e_i+e_j,\alpha_{A_2}\rangle=-1$, a contradiction 
because $e_i+e_j$ is in $H_{A_2}$. Thus $A_1\neq\emptyset$. Since all the 
vectors in ${\rm aff}({\cal A}')$ satisfy the linear 
equation
$$
\sum_{i=1}^m{x_i}=\sum_{i=m+1}^p{x_i},
$$
we obtain
$$
H_{A_2}^-\cap {\rm aff}({\cal A}')=
\left\{x\in{\rm aff}({\cal A}')\left\vert\, \ \ \sum_{v_i\in 
V_1\setminus N(A_2)}\hspace{-3mm}x_i\ \, \leq
\sum_{v_i\in V_2 \setminus A_2}\hspace{-1mm}x_i\right\}.\right.
$$
Hence we need only show $V_2\setminus A_2=N(A_1)$. The containment 
$N(A_1)\subset V_2\setminus A_2$ holds in general. For the reverse containment 
take $v_i\in V_2 \setminus A_2$. There is $v_j$ such 
that $\{v_i,v_j\}$ 
is an edge of $D$. If $v_j\in N(A_2)$, then 
$\langle e_i+e_j,\alpha_{A_2}\rangle=-1$, a contradiction 
because 
$e_i+e_j\in H_{A_2}$. Hence $v_j\notin N(A_2)$ 
and $v_i\in N(A_1)$. \QED

\medskip

For later use we state the following duality of 
facets which follows from the proof of Proposition~\ref{oct28-01}. 

\begin{Lemma}\label{dec28-01} Let $G$ be a 
connected bipartite graph with bipartition $(V_1,V_2)$ 
and let $F=H_A\cap\mathbb{R}_+{\cal A}$ be a facet of 
$\mathbb{R}_+{\cal A}$ with $A\subsetneq V_1$. Then 
\begin{description}
\item{\rm (a)} If $N(A)=V_2$, then $A=V_1\setminus\{v_i\}$ 
for some $v_i\in V_1$ and 
$F=H_{e_i}\cap\mathbb{R}_+{\cal A}$.\vspace{-1mm} 
\item{\rm (b)} If $N(A)\subsetneq V_2$, then 
$F=H_{V_2\setminus N(A)}\cap\mathbb{R}_+{\cal A}$ 
and $N(V_2\setminus N(A))=V_1\setminus A$.
\end{description}
\end{Lemma}

\begin{Definition}\rm For any set of vertices $S$ of 
a graph $G$, the {\it induced subgraph\/}
\index{induced subgraph} $\langle S\rangle$ is 
the maximal subgraph of $G$ with vertex set $S$.
\end{Definition}

\begin{Proposition}\label{dec27-01} Let $G$ be a 
connected bipartite graph with bipartition $(V_1,V_2)$ 
and let $A\subsetneq V_1$. Then 
$F=H_A\cap\mathbb{R}_+{\cal A}$ is a facet of 
$\mathbb{R}_+{\cal A}$ if and only if 
\begin{description}
\item{\rm (a)} $\langle A\cup N(A)\rangle$ is connected 
with vertex set $V(G)\setminus\{v\}$ for 
some $v\in V_1$, or\vspace{-1mm} 

\item{\rm (b)} $\langle A\cup N(A)\rangle$ and 
$\langle (V_2\setminus N(A))\cup (V_1\setminus A)\rangle$ are 
connected and their union is a spanning subgraph of $G$.
\end{description}
Moreover any facet has the form 
$F=H_A\cap\mathbb{R}_+{\cal A}$ for some $A\subsetneq V_i$, 
$i=1$ or $i=2$. 
\end{Proposition}

\demo The first statement follows readily 
from Lemma~\ref{dec28-01} 
and using part of the proof of Theorem~\ref{feasible1}. The last 
statement follows combining Theorem~\ref{feasible1} with 
Proposition~\ref{oct27-01}. \QED 

\begin{Remark}\rm In Proposition~\ref{dec27-01} 
the case (a) is included in case (b). To see this make 
$N(A)=V_2$ and note that 
$\langle (V_2\setminus N(A))\cup (V_1\setminus A)\rangle$ 
must consist of a point. The condition in case (a) 
is equivalent to require $G\setminus\{v\}$ 
connected and in this case 
$F=H_{e_i}\cap\mathbb{R}_+{\cal A}$, where $v=v_i$ 
correspond to the unit vector $e_i$.
\end{Remark} 

\begin{Lemma}\label{uniquf} Let $G$ be a connected 
bipartite graph with bipartition $(V_1,V_2)$ 
and let $F$ be a facet of $\mathbb{R}_+{\cal A}$. If 
$F=H_A\cap \mathbb{R}_+{\cal A}=H_B\cap
 \mathbb{R}_+{\cal A}$ with $A\subsetneq V_1$ and 
$B\subsetneq V_1$, then $A=B$. 
\end{Lemma}

\demo Set $V_1=\{v_1,\ldots,v_m\}$ and 
$V_2=\{v_{m+1},\ldots,v_{m+n}\}$. Recall that the equality
$$
x_1+\cdots+x_m=x_{m+1}+\cdots+x_{m+n}
$$
defines ${\rm aff}(\mathbb{R}_+{\cal A})$. 

Case (I): $N(A)=V_2$. Then by Lemma~\ref{dec28-01} (after permutation 
of vertices) $A=\{v_1,\ldots,v_{m-1}\}$. Hence any $x\in F$ satisfies
$$
\sum_{v_i\in A}x_i-\sum_{v_i\in N(A)}x_i=-x_m
$$
and thus $F=H_{e_m}\cap\mathbb{R}_+{\cal A}$. If $v_m\in B$, 
then $\{v_m,v_j\}\in E(G)$ for some $v_j$ in $N(B)$, 
thus $e_m+e_j\in H_B$ and consequently 
$e_m+e_j\in H_{e_m}$, a contradiction. Hence $v_m\notin B$, 
that is, $B\subset A$. If $N(B)=V_2$, then by Lemma~\ref{dec28-01} $A=B$. Assume 
$V_2\setminus N(B)\neq\emptyset$, to complete the proof for this 
case we will show that this assumption leads to 
a contradiction. First note that $v_m$ is not adjacent to 
any $v_j\in V_2\setminus N(B)$. Indeed if 
$\{v_m,v_j\}\in E(G)$, then $e_m+e_j\in H_B$. Thus 
$e_m+e_j\in H_{e_m}$, a contradiction. Therefore by the connectivity 
of $G$ at least one vertex $v_i\in V_1\setminus B$ must 
be adjacent to both a vertex $v_j\in V_2\setminus N(B)$ 
and a vertex $v_k\in N(B)$, which is impossible 
because $e_i+e_k\in H_{e_m}$ and $e_i+e_k\notin H_B$. 

Case (II): $N(A)\subsetneq V_2$ and $N(B)\subsetneq V_2$. 
We begin by considering the subcase $A\cap B\neq\emptyset$. 
Take $v_0\in A\cap B$ and $v_0\neq v\in B$. By Proposition~\ref{dec27-01} 
the subgraph $\langle B \cup N(B)\rangle$ is connected, hence 
there is a path of even length
$$
{\cal P}=\{v_0,v_1,v_2,\ldots,v_{2r-1},v_{2r}=v\}
$$
such that $v_{2i}\in B$ for all $i$. Note that $v_2\in A$. 
If $v_2\notin A$, then $e_1+e_2\in H_B$ and 
$e_1+e_2\notin H_A$, a contradiction. By induction we 
get $v_{2i}\in A\cap B$ for all $i$. Hence $v\in A$. This 
proves $B\subset A$, a similar argument proves $A=B$.  

Assume now $A\cap B=\emptyset$. We claim $N(A)\cap N(B)=\emptyset$, for otherwise 
if $\{v_j,v_k\}$ is an edge with $v_j\in B$ and 
$v_k\in N(A)\cap N(B)$, then $e_j+e_k\notin H_A$ because 
$v_j\notin A$ and $e_j+e_k\in H_B$, a contradiction. 

We may now assume 
$A\cap B=N(A)\cap N(B)=\emptyset$. Observe 
$A\cup B\neq V_2$ because if $V_2=A\cup B$, then $G$ 
would be disconnected with components 
$\langle A\cup N(A)\rangle$ and $\langle B\cup N(B)\rangle$. 
Take $v_j\in V_1\setminus(A\cup B)$ such that $v_j$ 
is adjacent to some $v_k$ in $N(A)\cup N(B)$, 
this choice is possible because $G$ is connected. Say 
$v_k\in N(A)$. Note $v_k\notin N(B)$. Then 
$e_j+e_k\in H_B$ and $e_j+e_k\notin H_A$, a 
contradiction. \QED

\medskip

Putting together the previous results we obtain 
the following canonical way of representing the edge cone.  The uniqueness 
follows from Lemma~\ref{dec28-01} and Lemma~\ref{uniquf}. 

\begin{Theorem}\label{uniqrepf} If $G$ is a connected 
bipartite graph with bipartition $(V_1,V_2)$, 
then there is a unique irreducible representation
$$
\mathbb{R}_+{\cal A}=
{\rm aff}(\mathbb{R}_+{\cal A})
\cap (\cap_{i=1}^rH_{A_i}^-)\cap(\cap_{i\in \cal I}
 H_{e_i}^+)  
$$
such that $A_i\subsetneq V_1$ for all $i$ and $v_i\in V_2$ for 
$i\in\cal{I}$. 
\end{Theorem}

\begin{Lemma}\label{oct14-01} If $G$ is a bipartite graph, then
$$
\mathbb{Z}^n\cap\mathbb{R}_+{\cal A}=\mathbb{N}{\cal A}.
$$
In particular if $(\beta_1,\ldots,\beta_n)$ is an integral vector 
in the edge cone,  then $\sum_{i=1}^n\beta_i$ is an even 
integer.
\end{Lemma}

\demo Let ${\cal A}=\{\alpha_1,\ldots,\alpha_q\}$ be the set of 
column vectors of the incidence matrix $M$ of $G$. Take 
$\alpha\in \mathbb{Z}^n\cap\mathbb{R}_+{\cal A}$, then by 
Carath\'eodory's Theorem \cite[Theorem~2.3]{ewald} and 
after an appropriate permutation of the $\alpha_i$'s we can write 
$$
\alpha=\eta_1\alpha_1+\cdots+\eta_r\alpha_r \ \ \ (\eta_i\geq 0),
$$
where $r$ is the rank of $M$ and $\alpha_1,\ldots,\alpha_r$ are 
linearly independent. On the other hand the submatrix 
$M'=(\alpha_1\cdots\alpha_r)$ is totally unimodular because $G$ is 
bipartite (see \cite{Schr}), hence by Kronecker's 
lemma \cite[p.~51]{Schr} the system of equations $M'x=\alpha$ has an 
integral solution. Hence $\alpha$ is a linear combination of 
$\alpha_1,\ldots,\alpha_r$ with coefficients in $\mathbb{Z}$. It follows 
that $\eta_i\in\mathbb{N}$ for all $i$, that is, 
$\alpha\in\mathbb{N}{\cal A}$. The other containment is clear. \QED

\medskip

As an application we recover the following version 
of the marriage problem for bipartite graphs, see \cite{Boll}. 
Recall 
that a pairing off of all the vertices of a graph $G$ is called 
a {\it perfect matching\/}\index{perfect matching\/}.

\begin{Theorem}[Marriage Theorem]\label{gen-marriage}
If $G$ is a bipartite graph, then  $G$ has a perfect matching if and only 
if 
$$|A|\leq |N(A)|$$ 
for every independent
set of vertices $A$ of $G$.\index{marriage theorem}
\end{Theorem}

\demo Note that $G$ has a perfect matching if and only if the 
vector $\beta=(1,1,\ldots,1)$ is in $\mathbb{N}{\cal A}$. By 
Lemma~\ref{oct14-01} $\beta$ is in $\mathbb{N}{\cal A}$ if 
and only if $\beta\in\mathbb{R}_+{\cal A}$. Thus the result 
follows from Corollary~\ref{feasible2}. \QED

\begin{Corollary}\label{kmn}
Let $G={\cal K}_{m,n}$ be the complete bipartite graph with $m\leq n$. If 
$V_1=\{v_1,\ldots,v_m\}$ and $V_2=V\setminus V_1$ is the bipartition 
of $G$, then
a vector $z\in {\mathbb R}^{m+n}$ is in ${\mathbb R}_+{\cal A}$ if 
and only if
$z=(x_1,\ldots,x_m,y_1,\ldots,y_n)$ satisfies  
\[
\begin{array}{rcll}
x_1+\cdots+x_m&=& y_1+\cdots+y_n, &
\\ -x_i& \leq & 0,& i=1,\ldots,m,\\
 -y_i& \leq & 0,& i=1,\ldots,n.
\end{array}
\]
In addition if $m\geq 2$, the inequalities define all the
facets of ${\mathbb R}_+{\cal A}$.
\end{Corollary}

\begin{Example}\rm 
If $G={\cal K}_{1,3}$ is the star with vertices $\{v,v_1,v_2,v_3\}$ and 
center $x$, then
the edge cone of $G$ has three facets defined by
\[
x_i\geq 0, \ \ \ (i=1,2,3).
\]
Note that $x=0$, define a proper face of dimension $1$.
\end{Example}

\end{document}